\documentclass[12pt,thmsa]{article}
\usepackage{amsfonts}
\usepackage{sw20lart}

\input tcilatex2.5
\begin{document}

\title{Geometric aspects of the Daugavet property.}
\author{R. V. Shvidkoy \\
Department of Mathematics\\
Mathematical Sciences Building\\
Columbia, Missouri, 65211\\
USA\\
\textit{e-mail: }mathgr31@showme.missouri.edu}
\date{November, 1998}
\maketitle

\begin{abstract}
Let $X$ be a closed subspace of a Banach space $Y$ and $J$ be the inclusion
map. We say that the pair $(X,Y)$ has the Daugavet property if for every
rank one bounded linear operator $T$ from $X$ to $Y$ the following equality 
\begin{equation}
\Vert J+T\Vert =1+\Vert T\Vert  \label{ade}
\end{equation}
holds. A new characterization of the Daugavet property in terms of weak open
sets is given. It is shown that the operators not fixing copies of $\ell _1$
on a Daugavet pair satisfy (\ref{ade}).

Some hereditary properties are found: if $X$ is a Daugavet space and $Y$ is
its subspace, then $Y$ is also a Daugavet space provided $X/Y$ has the
Radon-Nikod\'ym property; if $Y$ is reflexive, then $X/Y$ is a Daugavet
space. Becides, we prove that if $(X,Y)$ has the Daugavet property and $%
Y\subset Z$, then $Z$ can be renormed so that $(X,Z)$ possesses the Daugavet
property and the equivalent norm coincides with the original one on $Y$.
\end{abstract}

\section{Introduction.}

Let $X$ be a closed subspace of a Banach space $Y$ and $J:X\rightarrow Y$ be
the inclusion map. We say that the pair $(X,Y)$ has the Daugavet property
(or is a Daugavet pair) if for every rank one bounded linear operator $T$
from $X$ to $Y$ the following identity 
\begin{equation}
\Vert J+T\Vert =1+\Vert T\Vert ,  \label{DE}
\end{equation}
which is called the Daugavet equation, holds. If (\ref{DE}) is satisfied by
operators from some class $\mathcal{M}$ we say that $(X,Y)$ has the Daugavet
property with respect to this class.

(\ref{DE}) was first established for compact operators on $C[0,1]$ by
Daugavet in 1963 (see \cite{Daug}). Further it became a subject of extensive
study mostly directed to finding new Daugavet spaces and classes of
operators satisfying (\ref{DE}). In particular, it was proved that all
non-atomic $C(K)$ and $L_1(\mu )$-spaces possess the Daugavet property even
for weakly compact operators (see \cite{fs, h1, h2}). Until recently
investigation of general properties of Daugavet spaces remained somehow
aside. As far as we could trace the first results in this direction appeared
in works of Wojtaszczyk \cite{woj} and Kadets \cite{k}. Some of the most far
reaching ones were the following:

i) The unit sphere of a Daugavet space does not have a strongly exposed
point. Thus, a Daugavet space cannot have the Radon-Nikod\'ym property (see 
\cite{woj}).

ii) $\ell _1$ and $\ell _\infty $-sums of Daugavet spaces have the Daugavet
property (see \cite{woj} and \cite{kssw}).

iii) A Daugavet space does not have an unconditional basis (see \cite{k}).

A more intensive and systematic study of the general theory was initiated in 
\cite{kssw}. The authors gave a characterization of the Daugavet property in
terms of slices of the unit ball. This allowed to get a lot of information
about isomorphic structure of the Daugavet spaces.

The present paper is a natural continuation of \cite{kssw}. We give
affirmative answers for many questions posed there and provide alternative
proofs of some known earlier results.

In Section 2 another characterization of the Daugavet property in terms of
weak open sets intersecting the unit ball is given. Using this tool we prove
that all operators not fixing a copy of $\ell _1$ on a Daugavet pair satisfy
the Daugavet equation (Theorem \ref{oper}). Note that the analogous result
for strong Radon-Nikod\'ym operators was already obtained in \cite{kssw}. We
also present some new hereditary properties (Theorem \ref{her2}). In
particular, a pair $(X,Y)$ has the Daugavet property, provided $Y$ is a
Daugavet space and $Y/X$ has the Radon-Nikod\'ym property.

Section 3 is entirely devoted to pairs of the form $(X,C(K))$, where $K$ is
a compact Hausdorff space. It is shown that in some natural cases, e.g.,
when $K$ is the unit ball of $X^{*}$, such a pair possesses the Daugavet
property whenever $X$ does. We will see that this is also the case for some
bigger $C(K)$-spaces containing $X$. In Section 4 one of them is shown to
be, in a sense, universal: a Banach space $Y$ can be isomorphically embedded
into it, whenever $X\subset Y$ and $Y/X$ is separable.

At the end of Section 4 we prove following renorming theorem: let $(X,Y)$
have the Daugavet property and $Z$ be a Banach space containing $Y$, then $Z$
can be renormed so that $(X,Z)$ possesses the Daugavet property and the
equivalent norm remains unchanged on $Y$. A consequence of this result and
the aforementioned Theorem \ref{oper} is that a Daugavet space does not
embed into an unconditional sum of Banach spaces without copies of $\ell _1$%
. It is a generalization of the well known Theorem of Pelczy\'nski about
impossibility of embedding $C[0,1]$ and $L_1[0,1]$ into a space with
unconditional basis.

Throughout the text $\mathcal{L}(X,Y)$ denotes the space of all bounded
linear operators from $X$ into $Y$; $B(X)$ ($S(X)$) stands for the unit ball
(unit sphere) of a Banach space $X$; by $\overline{\limfunc{ext}}B(X^{*})$
we denote the weak$^{*}$ closure of the set of all extreme points of the
dual unit ball $B(X^{*}).$ For a subset $A$ of a Banach space, $\overline{A}$
denotes the norm-closure of $A$.

The author wishes to thank Professors V. Kadets, N. Kalton and D. Werner for
fruitful discussions, valuable remarks and constant interest to the work.

\section{Some characterizations and direct consequences.}

The central role in this section plays the notion of a slice.

\begin{definition}
Let $X$ be a Banach space. A $\emph{slice}$ of $B(X)$ is called the
following set 
\[
S(x^{*},\varepsilon )=\{x\in B(X):x^{*}(x)>1-\varepsilon \}, 
\]
where $x^{*}\in X^{*}$ and $\varepsilon >0$. We always assume that $x^{*}\in
S(X^{*})$. If $X$ is a dual space and $x^{*}$ is taken from the predual,
then $S(x^{*},\varepsilon )$ is called a weak$^{*}$ slice.
\end{definition}

In paper \cite{kssw} the following characterization of the Daugavet property
in terms of slices was obtained.

\begin{lemma}
\label{chl1}The following are equivalent:

\begin{enumerate}
\item[(a)]  The pair $(X,Y)$ has the Daugavet property;

\item[(b)]  For every $y_0\in S(Y)$ and for every slice $S(x_0^{*},%
\varepsilon _0)$ of $B(X)$ there is another slice $S(x_1^{*},\varepsilon
_1)\subset S(x_0^{*},\varepsilon _0)$ of $B(X)$ such that for every $x\in
S(x_1^{*},\varepsilon _1)$ the inequality $\Vert x+y_0\Vert \ge
2-\varepsilon _0$ holds;

\item[(c)]  For every $x_0^{*}\in S(X^{*})$ and for every weak$^{\,*}$ slice 
$S(y_0,\varepsilon _0)$ of $B(Y^{*})$ there is another weak$^{\,*}$ slice $%
S(y_1,\varepsilon _1)\subset S(y_0,\varepsilon _0)$ of $B(Y^{*})$ such that
for every $y^{*}\in S(y_1,\varepsilon _1)$ the inequality $\Vert
x_0^{*}+y_{\mid X}^{*}\Vert \ge 2-\varepsilon _0$ holds.
\end{enumerate}
\end{lemma}

For the sake of completeness we present the proof here.

\smallskip\emph{Proof.}

(a)$\Rightarrow $(b). Define $T{:}\allowbreak \ X\to Y$ by $Tx=x_0^{*}(x)y_0$%
. Then $\Vert J^{*}+T^{*}\Vert =\Vert J+T\Vert =2$, so there is a functional 
$y^{*}\in S(Y^{*})$ such that $\Vert J^{*}y^{*}+T^{*}y^{*}\Vert \ge
2-\varepsilon _0$ and $y^{*}(y_0)\ge 0$. Put 
\[
x_1^{*}=\frac{J^{*}y^{*}+T^{*}y^{*}}{\Vert J^{*}y^{*}+T^{*}y^{*}\Vert }%
,\quad \varepsilon _1=1-\frac{2-\varepsilon _0}{\Vert
J^{*}y^{*}+T^{*}y^{*}\Vert }. 
\]
Then for all $x\in S(x_1^{*},\varepsilon _1)$ we have 
\[
\langle (J^{*}+T^{*})y^{*},x\rangle \ge (1-\varepsilon _1)\Vert
J^{*}y^{*}+T^{*}y^{*}\Vert =2-\varepsilon _0, 
\]
therefore 
\begin{equation}
y^{*}(x)+y^{*}(y_0)x_0^{*}(x)\ge 2-\varepsilon _0,  \label{eq2}
\end{equation}
which implies that $x_0^{*}(x)\ge 1-\varepsilon _0$, i.e., $x\in
S(x_0^{*},\varepsilon _0)$. Moreover, by (\ref{eq2}) we have $%
y^{*}(x)+y^{*}(y_0)\ge 2-\varepsilon _0$ and hence $\Vert x+y_0\Vert \ge
2-\varepsilon _0$.

(b)$\Rightarrow $(a). Let $T\in \mathcal{L}(X,Y)$, $Tx=x_0^{*}(x)y_0$ be a
rank one operator. We can assume that $\Vert T\Vert =1$ (see, for example, 
\cite{AbraAB}) and $\Vert x_0^{*}\Vert =\Vert y_0\Vert =1$. Fix any $%
\varepsilon >0$. Then there is an $x\in S(x_0^{*},\frac \varepsilon 2)$ such
that $\Vert x+y_0\Vert >2-\frac \varepsilon 2$. So, 
\[
\Vert J+T\Vert \geq \Vert x+x_0^{*}(x)y_0\Vert \geq \Vert x+y_0\Vert
-|1-x_0^{*}(x)|>2-\varepsilon \text{.} 
\]
Let $\varepsilon $ go to zero.

The proof of equivalence (a)$\Leftrightarrow $(c) is analogous.\nopagebreak%
\hfill$\Box $

\bigskip One can see that the slices $S(x_1^{*},\varepsilon _1)$ and $%
S(y_1,\varepsilon _1)$ in the statement of Lemma \ref{chl1} can be replaced
by vectors $x$ and $y^{*}$. We will often refer to Lemma \ref{chl1} in this
form.

We mention some remarkable consequences of Lemma \ref{chl1} (the proofs can
be found in \cite{kssw}). First, if $X$ has the Daugavet property then $X$
(and $X^{*}$) contains an isomorphic copy of $\ell _1$, and moreover,
vectors equivalent to the canonical basis of $\ell _1$ can be chosen in
arbitrary slices of $B(X)$ (and weak$^{*}$ slices of $B(X^{*})$). Hence,
neither $X$ nor $X^{*}$ possess the Radon-Nikod\'ym property provided $X$
has the Daugavet property (see also \cite{woj} and \cite{wer}). Second, all
strong Radon-Nikod\'ym operators and, in particular, all weakly compact
operators on a Daugavet pair satisfy the Daugavet equation. Below we isolate
another such a class of operators, namely those not fixing copies of $\ell
_1 $, but first we need the following modification of Lemma \ref{chl1},
which shows that we can operate with weak open sets as well as with slices.

\begin{lemma}
\label{chl2}The following are equivalent:

\begin{enumerate}
\item[(a)]  The pair $(X,Y)$ has the Daugavet property;

\item[(b)]  For any given $\varepsilon >0$, $y\in S(Y)$ and weak open set $U$
in $X$ with $U\cap B(X)\neq \emptyset $ there is a weak open set $V$ in $X$
with $V\cap B(X)\neq \emptyset $ and $V\cap B(X)\subset U\cap B(X)$ such
that $\Vert v+y\Vert >2-\varepsilon $, whenever $v\in V\cap B(X)$;

\item[(c)]  For any given $\varepsilon >0$, $x^{*}\in S(X^{*})$ and weak$%
^{*} $ open set $U$ in $Y^{*}$ with $U\cap B(Y^{*})\neq \emptyset $ there is
a weak$^{*}$ open set $V$ in $Y^{*}$ with $V\cap B(Y^{*})\neq \emptyset $
and $V\cap B(Y^{*})\subset U\cap B(Y^{*})$ such that $\Vert v_{\mid
X}+x^{*}\Vert >2-\varepsilon $, whenever $v\in V\cap B(Y^{*})$.
\end{enumerate}
\end{lemma}

\emph{Proof.}\textbf{\ }Let us prove (a)$\Rightarrow $(b).

First we consider the weak$^{*}$ open set $U^{**}$ in $X^{**}$ that induces $%
U$ on $X$, i.e. $U^{**}\cap X=U$. By the Krein-Milman Theorem, there is a
convex combination of extreme points of $B(X^{**})$, $\stackrel{n}{%
\stackunder{i=1}{\sum }}\lambda _ix_i^{**}$, such that $\stackrel{n}{%
\stackunder{i=1}{\sum }}\lambda _ix_i^{**}\in U^{**}$. Clearly, we can find
weak$^{*}$ open neighborhoods $\{U_i^{**}\}_{i=1}^n$ of the points $%
\{x_i^{**}\}_{i=1}^n$ respectively, for which the following inclusion holds: 
\begin{equation}
\stackrel{n}{\stackunder{i=1}{\sum }}\lambda _i(U_i^{**}\cap
B(X^{**}))\subset U^{**}.  \label{incl1}
\end{equation}
Now by the Choquet Lemma (weak$^{*}$ slices containing an extreme point form
a basis of its weak$^{*}$ neighborhoods, \cite[p.49]{hhz}), we can assume
that the sets $\{U_i^{**}\cap B(X^{**})\}_{i=1}^n$ are weak$^{*}$ slices.
Thus, inclusion (\ref{incl1}) restricted on $X$ looks as follows: $\stackrel{%
n}{\stackunder{i=1}{\sum }}\lambda _iS_i\subset U$, where $S_i=U_i^{**}\cap
B(X^{**})\cap X$ are slices for all $i=1,2,\ldots ,n$.

Employing Lemma \ref{chl1}(b) we find a vector $x_1\in S_1$ with $\Vert
\lambda _1x_1+y\Vert >(\lambda _1+1-\varepsilon )$. Analogously, there is an 
$x_2\in S_2$ with $\Vert \lambda _2x_2+\lambda _1x_1+y\Vert >(\lambda
_2+\lambda _1+1-\varepsilon )$. Continuing in the same way we finally find $%
x_n\in S_n$ with $\Vert \lambda _nx_n+\lambda _{n-1}x_{n-1}+\ldots +\lambda
_1x_1+y\Vert >(\lambda _n+\lambda _{n-1}+\ldots +\lambda _1+1-\varepsilon
)=2-\varepsilon $, and $\stackrel{n}{\stackunder{i=1}{\sum }}\lambda
_ix_i\in U$. It remains only to use the lower weak semicontinuity of a norm
to get the required weak open set $V$.

This completes the proof of implication (a)$\Rightarrow $(b).

The implication (a)$\Leftarrow $(b) follows from Lemma \ref{chl1} and the
equivalence (a)$\Leftrightarrow $(c) is proved in the same way.\nopagebreak%
\hfill$\Box $

\begin{theorem}
\label{oper}If the pair $(X,Y)$ has the Daugavet property, then every
operator from $\mathcal{L}(X,Y)$ not fixing copies of $\ell _1$ satisfies
the Daugavet equation.
\end{theorem}

\emph{Proof.} Let $T\in \mathcal{L}(X,Y)$, $\Vert T\Vert =1$, be such an
operator and $\varepsilon >0$ be arbitrary.

Our considerations will rely on the following ``releasing principle'':
suppose for some finite set of vectors $\{x_i\}_{i=1}^n\subset B(X)$ and
some $\varepsilon >0$ the inequalities

\begin{equation}
\left\| \stackrel{n}{\stackunder{i=1}{\sum }}\theta _ix_i\right\|
>n-\varepsilon ,  \label{oper2}
\end{equation}
and 
\begin{equation}
\left\| \stackunder{i\in I_1}{\sum }a_ix_i+\stackunder{i\in I_2}{\sum }%
a_iTx_i\right\| >\left( \stackunder{i\in I_1\cup I_2}{\sum a_i}\right)
(1-\varepsilon )  \label{oper4}
\end{equation}
hold for all non-negative reals $a_i$, signs $\theta _i$, and some disjoint
sets $I_1$, $I_2\subset \{1,2,\ldots ,n\}$. Then there is a weak open set $%
U\subset X$ such that (\ref{oper2}) and (\ref{oper4}) remain true for all $%
x_n\in U\cap B(X)$.

Let us prove it. By the compactness argument, there is a $\delta >0$ such
that 
\begin{equation}
\left\| \stackunder{i\in I_1}{\sum }a_ix_i+\stackunder{i\in I_2}{\sum }%
a_iTx_i\right\| >1-\varepsilon +\delta ,  \label{oper3}
\end{equation}
whenever $\stackunder{i\in I_1\cup I_2}{\sum a_i}=1$ and $I_1$, $I_2$ as
above. Fix a finite $\frac \delta 2$-net $\left\{ (a_{k,1},a_{k,2},\ldots
,a_{k,n})\right\} _{k=1}^K$ in the set $\left\{ (a_1,a_2,\ldots ,a_n):~%
\stackrel{n}{\stackunder{i=1}{\sum }}a_i=1,~a_i\geq 0\right\} $ equipped
with the $\ell _1$-metric. Using the lower weak semicontinuity of a norm and
weak continuity of a bounded linear operator we conclude that there is a
weak open set $U$ such that both (\ref{oper2}) and (\ref{oper3}) hold for $%
a_i=a_{k,i}$ , $i=1,2,\ldots ,n$, $k=1,2,\ldots ,K$ and all $z_n\in U\cap
B(X)$. It is not hard to see that $U$ is desired.

Now we construct a sequence $\{x_i\}_{i=1}^\infty \subset B(X)$ which
satisfies (\ref{oper2}) and (\ref{oper4}) for all non-negative reals $a_i$,
signs $\theta _i$ and all disjoint finite sets $I_1$, $I_2\subset \Bbb{N}$ .

Assume that we have constructed such a sequence $\{x_i\}_{i=1}^n$ of length $%
n.$ We want to prove now that altering only the last term $x_n$ one can find
another vector $x_{n+1}$ such that the resulting sequence of length $n+1$
satisfies (\ref{oper2}) and (\ref{oper4}). Arguing in such a way, we produce
the desired infinite sequence if only take $x_1\in S(X)$ with $\Vert
Tx_1\Vert >1-\varepsilon $ on the first step.

Let us put $x_{n+1}^{\prime }=x_n$ for a moment. Clearly, (\ref{oper4})
remains true for the sequence $x_1,x_2,\ldots ,x_n,x_{n+1}^{\prime }$ and
all $I_1,I_2$ with additional restriction: if one of them contains $n$ then
the other does not contain $n+1$. We get rid of this restriction by
alteration of $x_n$ and $x_{n+1}^{\prime }$. To this end, we use the
`releasing principle' for $x_{n+1}^{\prime }$ and find the corresponding
weak open set $U\subset X$. Application of Lemma \ref{chl2}(b) several times
yields a vector $x_{n+1}\in U\cap B(X)$ such that (\ref{oper2}) is valid for
the sequence $x_1,x_2,\ldots ,x_n,x_{n+1}$ and (\ref{oper4}) holds without
the restriction: if $I_1$ contains $n+1$, then $I_2$ does not contain $n$.
Then we use the ``releasing principle'' to release $x_n$ so that both (\ref
{oper2}) and (\ref{oper4}) remain true. Appealing to Lemma \ref{chl2}(b) we
finally get an $x_n^{\prime }$ such that (\ref{oper4}) holds for the
sequence $x_1,x_2,\ldots ,x_n^{\prime },x_{n+1}$ without any restrictions on 
$I_1$ and $I_2$. Inequality (\ref{oper2}) is satisfied automatically.

The constructed sequence is $(1-\varepsilon )$-equivalent to the canonical
basis of $\ell _1$, for if $\stackrel{n}{\stackunder{i=1}{\sum }}|\lambda
_i|=1$, then by (\ref{oper2}) we have 
\begin{eqnarray*}
\left\| \stackrel{n}{\stackunder{i=1}{\sum }}\lambda _ix_i\right\|
&=&\left\| \stackrel{n}{\stackunder{i=1}{\sum }}\limfunc{sign}\lambda
_i\cdot x_i+\stackrel{n}{\stackunder{i=1}{\sum }}(\lambda _i-\limfunc{sign}%
\lambda _i)\cdot x_i\right\| \\
&>&n-\varepsilon -\stackrel{n}{\stackunder{i=1}{\sum }}|\lambda _i-\limfunc{%
sign}\lambda _i|=n-\varepsilon -\stackrel{n}{\stackunder{i=1}{\sum }}%
|1-|\lambda _i|| \\
&=&n-\varepsilon -n+1=1-\varepsilon .
\end{eqnarray*}

Since $T$ fixes no copies of $\ell _1$, by Rosenthal's Lemma we may assume
that the sequence $(Tx_n)_{n=1}^\infty $ is weakly Cauchy. Thus, $%
(Tx_{2n+1}-Tx_{2n})_{n=1}^\infty $ is weakly null. By Mazur's Theorem there
are two finite disjoint sets $I_1$, $I_2\subset N$ such that for some $p\in 
\limfunc{conv}\{x_i:i\in I_1\}$ and $q\in \limfunc{conv}\{x_i:i\in I_2\}$ we
have $\Vert Tp-Tq\Vert <\varepsilon $. From this and (\ref{oper4}) we
finally obtain 
\[
\Vert p+Tp\Vert >\Vert p+Tq\Vert -\varepsilon >2(1-\varepsilon )-\varepsilon
=2-3\varepsilon , 
\]
which implies $\Vert J+T\Vert =2$ in view of arbitrariness of $\varepsilon $.

This finishes the proof.\nopagebreak\hfill$\Box $

\bigskip
It is known that $C(K)$ has the Daugavet property (see \cite{fs} or \cite{h2}%
) if $K$ is a compact Hausdorff space without isolated points. Besides, due
to a result of Rosenthal \cite{ros} and by Lemma 2.4 from \cite{ww} it
follows that operators on $C(K)$ not fixing copies of $C[0,1]$ are precisely
those not fixing copies of $\ell _1$. So, from the previous theorem we
obtain that all such operators satisfy the Daugavet equation. This result
was first established by Weis and Werner in their paper \cite{ww}. By
Theorem \ref{oper} we also solve a problem posed in \cite{kssw}.

\begin{corollary}
Suppose $X$ is a Daugavet space and $Y$ is a complemented subspace in $X$
such that $X/Y$ contains no copies of $\ell _1$, then the norm of every
projection from $X$ onto $Y$ is at least $2$.
\end{corollary}

\emph{Proof.} Let $P:X\rightarrow X$ be any projection onto $Y$. Then $-Id+P$
fixes no copies of $\ell _1$ and hence, by Theorem \ref{oper}, satisfies the
Daugavet equation. So, we have $\Vert P\Vert =\Vert Id+(-Id+P)\Vert =1+\Vert
P-Id\Vert \geq 2$.\nopagebreak\hfill$\Box $

\bigskip\noindent \textbf{Problem 1. }It remains open whether every
Dunford-Pettis operator on a Daugavet pair satisfies the Daugavet equation.

\smallskip\noindent \textbf{Problem 2. }One of the remarkable
characterizations of Banach spaces not containing isomorphic copies of $\ell
_1$ is that the duals of such spaces possess the weak Radon-Nikod\'ym
property. Thus, no dual to a Daugavet space has this property. It is not
known, however, if the same is true for a Daugavet space itself.

\bigskip Now we discuss the following question: suppose $X$ has the Daugavet
property; what classes of subspaces of $X$ possess the same property?

It was shown in\textrm{\ }\cite{kssw} that all the subspaces with separable
annihilator do. Such an effect could be attributed to extreme ``spreadness''
of a Daugavet unit ball (see Lemmas \ref{chl1} and \ref{chl2}). We will
repeatedly use this idea later on.

\begin{theorem}
\label{her2}Let $X$ have the Daugavet property and $Y$ be a subspace of $X$.

\begin{enumerate}
\item[(a)]  If $X/Y$ has the Radon-Nikod\'ym property, then the pair $(Y,X)$
has the Daugavet property;

\item[(b)]  If $Y$ is reflexive, then $X/Y$ has the Daugavet property.
\end{enumerate}
\end{theorem}

In the particular case when $X=L_1[0,1]$ part (b) of~Theorem \ref{her2} was
proved in \cite{kssw}.

\smallskip\emph{Proof.} Part (a). According to Lemma \ref{chl1}(b) it is
sufficient to prove that given any $\delta >0$, $S(y^{*},\varepsilon )$ and $%
x\in B_X$ ~there is a $y\in S(y^{*},\varepsilon )$ such that $\Vert x+y\Vert
>2-\delta .$

Denote by $j$ the quotient map $:X\mapsto X/Y$. Saving the notation for the
functional $y^{*}$, we extend it to all of $X$ by the Hahn-Banach Theorem.
The set $A=j(S(y^{*},\varepsilon ))$ is convex and contains the origin.
Since $X/Y$ has the Radon-Nikod\'ym property, the Phelps Theorem (see for
example \cite{Die-LNM}) yields a convex combination $\stackrel{n}{%
\stackunder{i=1}{\sum }}\lambda _ia_i$ of strongly exposed points $%
\{a_i\}_{i=1}^n$ of the set $\overline{A}$ for which 
\begin{equation}
\left\| \stackrel{n}{\stackunder{i=1}{\sum }}\lambda _ia_i\right\| <\frac
\delta 2.  \label{hereq1}
\end{equation}

Let $\{a_i^{*}\}_{i=1}^n\subset (X/Y)^{*}$ be functionals exposing $%
\{a_i\}_{i=1}^n$ respectively and let positive numbers $\{\varepsilon
_i\}_{i=1}^n$ be such that 
\begin{equation}
\limfunc{diam}\left\{ S(a_i^{*},\varepsilon _i)\cap \overline{A}\right\}
<\frac \delta 4\text{, }i=1,2,...,n.  \label{hereq2}
\end{equation}
Since $S(a_i^{*},\varepsilon _i)\cap A\neq \emptyset $, we have $%
S(j^{*}a_i^{*},\varepsilon _i)\cap S(y^{*},\varepsilon )\neq \emptyset $.
Applying Lemma \ref{chl2}(b) we find $x_i\in S(j^{*}a_i^{*},\varepsilon
_i)\cap S(y^{*},\varepsilon )$ such that 
\[
\left\| \stackrel{n}{\stackunder{i=1}{\sum }}\lambda _ix_i+x\right\|
>2-\frac \delta 4 
\]
Now taking into account (\ref{hereq1}) and (\ref{hereq2}) we obtain the
following estimate: 
\[
\left\| j(\stackrel{n}{\stackunder{i=1}{\sum }}\lambda _ix_i)\right\|
<\left\| \stackrel{n}{\stackunder{i=1}{\sum }}\lambda _ia_i\right\| +\frac
\delta 4<\frac \delta 2. 
\]
It means that there is a $y\in B_Y$ for which 
\[
\left\| \stackrel{n}{\stackunder{i=1}{\sum }}\lambda _ix_i-y\right\| <\delta
. 
\]
Then by (\ref{hereq2}) we finally get 
\[
\Vert x+y\Vert >2-\frac 32\delta . 
\]
Clearly, $y\in S(y^{*},\varepsilon +\delta )$.

Because of arbitrariness of $\varepsilon $ and $\delta $, part (a) is proved.

The proof of part (b) is analogous (we have only to use the weak$^{*}$
topology and apply Lemma \ref{chl2}(c)). \nopagebreak\hfill$\Box $

\bigskip\noindent \textbf{Problem 3.} Under the conditions of Theorem \ref
{her2},

(a) does $Y$ have the Daugavet property if $X/Y$ is an Asplund space
(equivalently, $(X/Y)^{*}$ has the Radon-Nikod\'ym property) or, more
generally, if $X/Y$ fails to contain isomorphic copies of $\ell _1$?

(b) does $X/Y$ have the Daugavet property if either $Y$ or $Y^{*}$ (or both)
has the Radon-Nikod\'ym property or fails to contain isomorphic copies of $%
\ell _1$?

\section{Subspaces of C(K)-spaces.}

Now we study the case when in a pair $(X,Y)$ the space $Y$ is a $C(K)$-space
for some compact Hausdorff space $K$. As was shown in various works (see 
\cite{fs} or \cite{ww}) and as also follows from our Lemma \ref{chl1}, $C(K)$
has the Daugavet property if and only if $K$ has no isolated points.
Moreover, we can assert that if for some $X\subset C(K)$ the pair $(X,C(K))$
has the Daugavet property, then $K$ does not have such a point $k$, for
otherwise the rank one operator $Tx=-\chi _{\{k\}}\cdot x(k)$ does not
satisfy the Daugavet equation. So, investigating pairs of the form $(X,C(K))$
it is natural to require that $K$ have no isolated points.

We begin with a characterization of those Banach spaces $X$, $X\subset C(K)$%
~that the pair $(X,C(K))$ has the Daugavet property. In the sequel, $\delta
_k^{*}$, $k\in K$ stands for the functional on $C(K)$ acting by the rule $%
\delta _k^{*}(f)=f(k)$, $f\in C(K)$.

\begin{lemma}
\label{chC(K)}Let $X$ be a subspace of $C(K)$, where $K$ is a compact
Hausdorff space without isolated points. The following conditions are
equivalent:

\begin{enumerate}
\item[(a)]  The pair $(X,C(K))$ has the Daugavet property;

\item[(b)]  For every $\varepsilon >0$, $x^{*}\in S(X^{*})$ and open set $U$
in $K$ there exists a point $u\in U$ such that $\Vert x^{*}+\delta _{u\mid
X}^{*}\Vert >2-\varepsilon $;

\item[(c)]  For every $x^{*}\in S(X^{*})$ and open set $U$ in $K$ there
exists a (closed) $G_\delta $-set $G$ in $U$ such that $\Vert x^{*}+\delta
_{u\mid X}^{*}\Vert =2$, whenever $u\in G$.
\end{enumerate}
\end{lemma}

\emph{Proof.} (a)$\Rightarrow $(b). Let $f\in S(C(K))$ be a function
vanishing outside $U$. By Lemma \ref{chl1}(c), there is a slice $S\subset
S(f,\frac 12)$ such that $\Vert x^{*}+\mu \Vert >2-\varepsilon $, for all $%
\mu \in S$. Pick any $\delta _u^{*}\in S$. Clearly, $\delta
_u^{*}(f)=f(u)>\frac 12$ and hence, $u\in U$. So, $u$ is the required point.

(b)$\Rightarrow $(c). Apply part (b) countably many times and use the weak$%
^{*}$ lower semicontinuity of a dual norm and the regularity of a Hausdorff
compact set.

(c)$\Rightarrow $(a). We apply Lemma \ref{chl1} again. Pick arbitrary $%
x^{*}\in S(X^{*})$ and weak$^{*}$ slice $S(f,\varepsilon )$ in $B(C^{*}(K))$%
. Let $U=\{k\in K:f(k)>1-\varepsilon \}$. By condition (c), we can find a
point $u\in U$ such that $\Vert x^{*}+\delta _{u\mid X}^{*}\Vert =2$.
Moreover, we have $\delta _u^{*}(f)=f(u)>1-\varepsilon $ and hence, $\delta
_u^{*}\in S(f,\varepsilon ).$ This completes the proof.\nopagebreak\hfill$%
\Box $

\bigskip Of course, not every pair $(X,C(K))$ has the Daugavet property
provided $X$ does, e.g., this one $(C[0,1],C([0,1]\cup [2,3]))$. However, as
the following theorem shows, in some natural and useful cases this is true.

\begin{proposition}
\label{prC(K)1}If the pair $(X,Y)$ has the Daugavet property and $K$ is
either $B(Y^{*})$ or $\overline{\limfunc{ext}}B(Y^{*})$, then the pair $%
(X,C(K))$ also has the Daugavet property.
\end{proposition}

\emph{Proof.} In both cases we use condition (b) of Lemma \ref{chC(K)}.

First, consider $K=B(Y^{*})$. Fix arbitrary $\varepsilon >0$, open set $%
U\subset K$ and $x^{*}\in S(X^{*})$. By Lemma \ref{chl2}(c) there is $%
y^{*}\in U$ such that $\Vert x^{*}+y_{\mid X}^{*}\Vert >2-\varepsilon $. We
denote by $u$ the functional $y^{*}$ regarding it as a point of topological
space $K$. It remains to notice that $\delta _{u\mid X}^{*}=y_{\mid X}^{*}$.

Let $K=\overline{\limfunc{ext}}B(X^{*})$. Fix $\varepsilon $, $U$ and $x^{*}$
as above. By the Choquet Lemma we may assume that $U$ is induced by a slice $%
S$. By Lemma \ref{chl1}(c) there is a slice $S_1\subset S$, and hence, there
is a $y^{*}\in S\cap K$ such that $\Vert x^{*}+y_{\mid X}^{*}\Vert
>2-\varepsilon $. So, as above the point $u=y^{*}$ is required.\nopagebreak%
\hfill$\Box $

\bigskip In the case $K=B(Y^{*})$ this proposition solves a problem posed in 
\cite{kssw}. The result was proved there for $K=\overline{\limfunc{ext}}%
B(Y^{*})$. However, we include both cases to emphasize their common origin.

\smallskip Let $K$ be a compact Hausdorff space without isolated points. We
introduce the following spaces: 
\begin{eqnarray*}
l_\infty (K) &=&\left\{ f:K\mapsto R,\quad \Vert f\Vert _\infty =\sup
(|f(s)|,\,s\in K)<{\infty }\right\} , \\
m(K) &=&\left\{ f\in l_\infty (K):\limfunc{supp}(f)\ \text{is a first
category set}\right\} , \\
m_0(K) &=&l_\infty (K)/m(K).
\end{eqnarray*}

In what follows we investigate Daugavet properties of the space $m_0(K)$. In
the next section we use them to prove some general results on renormings.

$m_0(K)$ equipped with the factor-norm is a real $C^{*}$-algebra, and hence,
is a $C(Q)$-space. The appropriate compact set $Q=Q_K$ can be defined as the
set of all real homomorphisms on $m_0(K)$ endowed with the induced weak$^{*}$
topology. This is precisely limits by ultrafilters on $K$, which do not
contain first category sets. Let $\frak{U}$ be such an ultrafilter. We
denote by $\lim \frak{U}$ the point in $K$ to which it converges and by $%
\frak{U}\_\lim $ the real homomorphism on $m_0(K)$ it generates ($\frak{U}%
\_\lim \in Q_K$).

\begin{lemma}
\label{al}Suppose $U$ is an open set in $Q_K$, then there is an open set $V$
in $K$ such that for every $v\in V$ one can find an ultrafilter $\frak{U}_v$
on $K$ with $\lim \frak{U}_v=v$ and $\frak{U}_v\_\lim \in U$.
\end{lemma}

\emph{Proof.} By the construction of $Q_K$ we may assume there are a finite
set $(f_i)_{i=1}^n\subset m_0(K)$, $\varepsilon >0$ and ultrafilter $\frak{U}%
_0$ on $K$ such that $U=\{\varphi \in Q_K:|\varphi (f_i)-\frak{U}_0\_\lim
(f_i)|<\varepsilon \}$. Denote $a_i=\frak{U}_0\_\lim (f_i)$. We fix a second
category set $A\in \frak{U}_0$ with the following property: 
\begin{equation}
f_i(A)\subset (a_i-\varepsilon \;,\;a_i+\varepsilon )\text{,\quad }%
i=1,2,\ldots ,n.  \label{incl}
\end{equation}
Then we find an open set $V$ in $K$ such that for any open $W\subset V$, $%
W\cap A$ is a second category set (see \cite{Kelley}). It remains to show
that $V$ is required.

Indeed, let $v\in V$. Consider an ultrafilter $\frak{U}_v$ containing $%
\{W\cap A:W$ is an open neighborhood of $v\}$. Plainly, $\lim \frak{U}_v=v$.
On the other hand, in view of (\ref{incl}) we have $\frak{U}_v\_\lim
(f_i)\in (a_i-\varepsilon \;,\;a_i+\varepsilon )$, $i=1,2,\ldots ,n$. This
means that $\frak{U}_v\_\lim \in U$. This finishes the proof.\nopagebreak%
\hfill$\Box $

\smallskip It is easy to see that $C(K)$ is isometrically embedded into $%
m_0(K)$ by the quotient map.

\begin{proposition}
\label{prC(K)2}If the pair $(X,C(K))$ has the Daugavet property, then the
pair $(X,m_0(K))$ also has the Daugavet property.
\end{proposition}

\emph{Proof.} We apply Lemma \ref{chC(K)} again using the interpretation of $%
m_0(K)$ as a $C(Q)$-space. To this end, we fix $\varepsilon >0$, open set $%
U\subset Q_K$ and $x^{*}\in S(X^{*})$. Applying Lemma \ref{al} to $U$ we
find the corresponding open set $V\subset K$. Lemma \ref{chC(K)} applied to
the pair $(X,C(K))$ yields $v\in V$ such that $\Vert x^{*}+\delta _{v\mid
X}^{*}\Vert >2-\varepsilon $. Consider the ultrafilter $\frak{U}_v$ with $%
\lim \frak{U}_v=v$ and $\frak{U}_v\_\lim \in U$, and denote $u=\frak{U}%
_v\_\lim $. So, $\delta _{v\mid X}^{*}=\delta _{u\mid X}^{*}$ and $u\in U$.
Hence, the point $u$ is desired.\nopagebreak\hfill$\Box $

\begin{corollary}
The pair $(C(K),m_0(K))$ has the Daugavet property.\hfill$\Box $
\end{corollary}

\begin{corollary}
\label{corm(K)}Let the pair $(X,Y)$ have the Daugavet property and $K$ be
either $B(Y^{*})$ or $\overline{\limfunc{ext}}B(X^{*})$, then the pair $%
(X,m_0(K))$ has the Daugavet property too.
\end{corollary}

\emph{Proof.} Combine Propositions \ref{prC(K)1} and \ref{prC(K)2}.%
\nopagebreak\hfill$\Box $

\section{Renorming theorem.}

The main goal of this section is to prove the following result.

\begin{theorem}
\label{renth}Let $X$, $Y$, $Z$ be Banach spaces such that $X\subset Y\subset
Z$. If the pair $(X,Y)$ has the Daugavet property, then $Z$ can be renormed
so that $(X,Z)$ possesses the Daugavet property and the equivalent norm
coincides with the original one on $Y$.
\end{theorem}

In separable case this theorem was proved in \cite{kssw}. The general case,
however, requires more detailed consideration. Therefore we present the
complete proof here.

First we prove a theorem which establishes, in some sense, a property of
universality of $m_0(K)$-spaces, where $K$ is the unit ball of a dual space.
Since in the sequel we often deal with density character of a Banach space $%
X $ (the minimal cardinality of a dense set in $X$), we denote it by $%
\limfunc{dens}(X)$.

\begin{theorem}
\label{embth}Let $Y$ be a closed subspace of Banach spaces $Z$ and $W$. Let
also $\limfunc{dens}(Z/Y)=\beta $, where $\beta $ is an ordinal. Suppose $%
B(W^{*})$ contains a family $\{B_\alpha \}_{\alpha <\beta }$ of disjoint
second category sets such that if $B^{\prime }=\stackunder{\alpha <\beta }{%
\cup }B_\alpha $, then $B^{\prime }\cap -B^{\prime }=\emptyset $. Then there
is an isomorphic embedding $E:Z\rightarrow m_0(B(W^{*}))$, which coincides
with the natural one on $Y$.
\end{theorem}

\emph{Proof.} Let us fix a dense set $([z_\alpha ])_{\alpha <\beta }\subset
B(Z/W)$ with $\Vert z_\alpha \Vert \leq 1$, and for every $\alpha <\beta $
find a functional $\varphi _\alpha \in S(Y^{\bot })$ so that $\varphi
_\alpha (z_\alpha )=\Vert [z_\alpha ]\Vert $. Also to every $w^{*}$ we
assign a functional $\widetilde{w}^{*}$ obtained by restriction of $w^{*}$
on $Y$ and then extension to all of $Z$ by the Hahn-Banach Theorem.

Now we want to embed $Z$ into $\ell _\infty (B(W^{*}))$ so that every
element from the image of $B(Z)$ takes values greater than $\frac 18$ on a
second category set. To this end, for each $z\in Z$ we define a function $%
f_z\in \ell _\infty (B(W^{*}))$ as follows: 
\[
f_z(w^{*})=\left\{ 
\begin{array}{ll}
\widetilde{w}^{*}(z), & w^{*}\in B(W^{*})\backslash B_0 \\ 
\widetilde{w}^{*}(z)+8\varphi _\alpha (z), & w^{*}\in B_\alpha
\end{array}
\right. . 
\]
Clearly the mapping $F:z\rightarrow f_z$ is linear and bounded. Moreover, $%
f_z(w^{*})=w^{*}(z)$, if $z\in Y$. So, $F_{\mid X}$ is the natural embedding
of $Y$ into $\ell _\infty (B(W^{*}))$ (even into $C(B(W^{*}))$).

Suppose now $\Vert z\Vert =1$. Then either $\Vert [z]\Vert \leq \frac 14$ or 
$\Vert [z]\Vert >\frac 14$. In the former case there is a $y_0\in Y$ such
that $\Vert z-y_0\Vert <\frac 38$. Because of the condition imposed on $%
B^{\prime }$, the set $\{w^{*}\in B(W^{*})\backslash B^{\prime
}:w^{*}(y_0)>\Vert y_0\Vert -\frac 18\}$ is of second category, and for
every its element we have 
\[
|f_z(w^{*})|=|\widetilde{w}^{*}(z)|>|\widetilde{w}^{*}(y_0)|-\frac
38=|w^{*}(y_0)|-\frac 38=\Vert y_0\Vert -\frac 12>\frac 18. 
\]
So, $|f_z(w^{*})|>\frac 18$, for $w^{*}$ from some second category set.

In the case $\Vert [z]\Vert >\frac 14$, there is an ordinal $\alpha $, $%
\alpha <\beta $, and $y\in Y$ such that $\Vert [z_\alpha ]\Vert >\frac 14$
and $\Vert z-z_\alpha -y\Vert <\frac 1{16}$. From this we get for all $%
w^{*}\in B_\alpha $%
\begin{eqnarray*}
|f_z(w^{*})| &=&|\widetilde{w}^{*}(z)+8\varphi _\alpha ^{*}(z)| \\
&>&|8\varphi _\alpha ^{*}(z_\alpha -y)|-\frac 12-|\widetilde{w}%
^{*}(z)|=8\Vert [z_\alpha ]\Vert -\frac 32 \\
&>&\frac 84-\frac 32=\frac 12.
\end{eqnarray*}

To define the desired isomorphic embedding $E:Z\rightarrow m_0(B(W^{*}))$ we
just put $Ez=[Fz]$, $z\in Z$.\nopagebreak\hfill$\Box $

It is not hard to construct countable number of second category sets
satisfying the condition of the previous theorem. So, in the special case
when $Z/Y$ is separable, we obtain the following corollary.

\begin{corollary}
Let $Y$ be a closed subspace of $Z$ such that $Z/Y$ is separable. Then there
exists an isomorphic embedding of $Z$ into $m_0(B(Y^{*}))$, which coincides
with the natural one on $Y$.
\end{corollary}

\emph{Proof of Theorem \ref{renth}.}

Suppose $(X,Y)$ is a Daugavet pair and $Z$ is some Banach space containing $%
Y $. If $B(Y^{*})$ were very ``reach'' of disjoint second category sets,
i.e. enough to satisfy the condition of Theorem \ref{embth} (in this case $%
Y=W$), there would exist an isomorphic embedding $E$ of $Z$ into $%
m_0(B(Y^{*}))$. Appealing to Corollary \ref{corm(K)}, the equivalent norm $%
|||z|||=\Vert Ez\Vert $ would be desired.

That, however, may not be the case, for example, when $\limfunc{dens}(Z)>%
\limfunc{dens}(m_0(B(Y^{*})))$. So, we should replace $Y$ by a bigger space,
say $W$, which meets the condition of Theorem \ref{embth} and at the same
time possesses the Daugavet property in pair with $X$. If we can do this,
the norm introduced in the previous case satisfies our requirements, and we
are done.

Let $\beta $ be as in Theorem \ref{embth}. We define $W$ to be the $\ell
_\infty $-sum of $\beta $ copies of $C(B(Y^{*}))$, i.e. $W=\left\{ (f_\alpha
)_{\alpha <\beta }:f_\alpha \in C(B(Y^{*}))\text{ and }\Vert (f_\alpha
)\Vert =\stackunder{\alpha <\beta }{\sup }\Vert f_\alpha \Vert <\infty
\right\} $. $Y$ embeds into $W$ as follows:

\begin{eqnarray*}
y &\rightarrow &(y_\alpha )_{\alpha <\beta }\text{,\quad }y\in Y\text{;} \\
y_\alpha (s) &=&s(y)\text{,\quad }s\in B(Y^{*})\text{.}
\end{eqnarray*}
So, $Y$ can be regarded as a subspace of $W$. Using Proposition \ref{prC(K)1}%
, it is not difficult to prove that the pair $(X,W)$ has the Daugavet
property.

Now fix $f\in C(B(Y^{*}))$, $\Vert f\Vert =1$, and for every $\alpha $, $%
\alpha <\beta $, define the vector $w_\alpha =(f_{\alpha ^{\prime
}})_{\alpha ^{\prime }<\beta }$ so that $f_{\alpha ^{\prime }}=f$, if $%
\alpha ^{\prime }=\alpha $, and $f_{\alpha ^{\prime }}=0$ otherwise. Put $%
B_\alpha =S(w_\alpha ,\frac 13)$. Since every $B_\alpha $ is weak$^{*}$
open, it is a second category set. Next, $B_{\alpha ^{\prime }}\cap
B_{\alpha ^{\prime \prime }}=\emptyset $, $\alpha ^{\prime }\neq \alpha
^{\prime \prime }$, for otherwise every $w^{*}\in B_{\alpha ^{\prime }}\cap
B_{\alpha ^{\prime \prime }}$ would have norm bigger than 1. For the same
reason, $B^{\prime }=\stackunder{\alpha <\beta }{\cup }B_\alpha $ is
disjoint with $-B^{\prime }$.

So, we have constructed the space satisfying all our requirements. This
finishes the proof. \nopagebreak\hfill$\Box $

\begin{corollary}
\label{corunc}A Daugavet space does not isomorphically embed into an
unconditional sum of Banach spaces without copies of $\ell _1$.
\end{corollary}

The proof is the same as that of Corollary 2.7 in \cite{kssw}. We only have
to use our Theorem \ref{oper} and the fact that the sum of finite number
operators not fixing copies of $\ell _1$ is an operator not fixing copies of 
$\ell _1$.

\smallskip It is worthwhile to remark that the previous result is a direct
generalization of the known Theorem of Pelczy\'nski for $C[0,1]$ and $%
L_1[0,1]$ spaces (for more about that see \cite{En-St}, \cite{ks} and \cite
{Kalton}).

\bigskip\noindent \textbf{Problem 4.} It would be interesting to find answer
to the following question: if $(X,Y)$ is a Daugavet pair, can $Y$ be
renormed to have the Daugavet property. We may, however, assert that such a
renorming cannot be accomplish leaving the norm on $X$ unchanged. In fact,
look at the space $L_\infty [0,1]$. It is 1-complemented in every containing
Banach space. Since every 1-codimensional subspace of a Daugavet space is at
least 2-complemented, $L_\infty [0,1]\oplus $ $\Bbb{R\ }$cannot be renormed
to have the Daugavet property so that the equivalent norm remains the same
on $L_\infty [0,1]$.


\begin{thebibliography}{99}
\bibitem{AbraAB}  Y.~Abramovich, C.~D. Aliprantis, and O.~Burkinshaw, %
\newblock The Daugavet equation in uniformly convex Banach spaces, 
\newblock 
\textit{J. Funct. Anal.} \newblock \textbf{97} (1991), 215--230.

\bibitem{Daug}  I.~K. Daugavet, \newblock On a property of completely
continuous operators in the space $C$, \newblock \textit{Uspekhi Mat. Nauk }%
\newblock \textbf{18.5} (1963), 157--158 (Russian).

\bibitem{Die-LNM}  J.~Diestel, \newblock ``Geometry of Banach Spaces --
Selected Topics,'' \newblock Lecture Notes in Math. 485. Springer,
Berlin-Heidelberg-New York, 1975.

\bibitem{En-St}  P. Enflo and T. W. Starbird, \newblock  Subspaces of $L_1$
containing $L_1$, \newblock  \textit{Studia Math. }\newblock  \textbf{65}
(1979), 203-225.

\bibitem{fs}  C.~Foia\c s and I.~Singer, \newblock Points of diffusion of
linear operators and almost diffuse operators in spaces of continuous
functions, \newblock \textit{Math. Z.} \newblock \textbf{87} (1965),
434--450.

\bibitem{hhz}  P.~Habala, P.~H\'ajek, and V.~Zizler, \newblock
``Introduction to Banach Spaces,'' \newblock Matfyz Press, Prague, 1996.

\bibitem{h1}  J.~R. Holub, \newblock Daugavet's equation and operators on $%
L_1(\mu )$, \newblock \textit{Proc. Amer. Math. Soc.} \newblock
\textbf{100} (1987), 295--300.

\bibitem{h2}  J.~R. Holub, \newblock A property of weakly compact operators
on $C[0,1]$,\newblock \textit{Proc. Amer. Math. Soc.} \newblock
\textbf{97} (1986), 396--398.

\bibitem{k}  V.~M. Kadets, \newblock Some remarks concerning the Daugavet
equation, \newblock \textit{Quaestiones Math.} \newblock
\textbf{19} (1996), 225--235.

\bibitem{kp}  V.~M. Kadets and M.~M. Popov, \newblock The Daugavet property
for narrow operators in rich subspaces of $C[0,1]$ and $L$\emph{$_1[0,1]$}, %
\newblock \textit{St. Petersburg Math. J.} \newblock \textbf{8} (1996),
43--62.

\bibitem{ks}  V.~M. Kadets and R.~V. Shvidkoy, \newblock The Daugavet
property for pairs of Banach spaces, \newblock \textit{Math. Analysis,
Algebra and Geometry} \newblock  (to appear).

\bibitem{kssw}  V.~M. Kadets, R.~V. Shvidkoy, G. G. Sirotkin, and D. Werner, %
\newblock Banach spaces with the Daugavet property, \newblock \textit{%
C.R.Acad.Sci.Paris} \newblock  \textbf{325.1} (1997), 1291-1294.

\bibitem{Kalton}  N. J. Kalton, \newblock The endomorphisms of $L_p$ ($0\leq
p\leq 1$), \newblock  \textit{Indiana Univ. Math. J.} \newblock  Vol. 27, 
\textbf{3 (}1978), 353-381.

\bibitem{Kelley}  J.~L. Kelley, \newblock ``General Topology,'' \newblock %
Van Nostrand, 1955.

\bibitem{ros}  H. P. Rosenthal, \newblock On factors of $C[0,1]$\ with
nonseparable dual, \newblock \textit{Israel J. of Math.} \newblock \textbf{13%
} (1972), 361-378.

\bibitem{ww}  L.~Weis and D.~Werner, \newblock The Daugavet equation for
operators not fixing a copy of $C[0,1]$, \newblock \textit{J. Operator Theory%
} \newblock  \textbf{39 }(1998), 89-98.

\bibitem{wer}  D.~Werner, \newblock The Daugavet equation for operators on
function spaces, \newblock \textit{J. Funct. Anal.} \newblock
\textbf{143} (1997), 117--128.

\bibitem{woj}  P.~Wojtaszczyk, \newblock Some remarks on the Daugavet
equation, \newblock \textit{Proc. Amer. Math. Soc.} \newblock \textbf{115}
(1992), 1047--1052.
\end{thebibliography}
\end{document}